# Asymptotic evaluation of a function defined by power series

Mihail M. Nikitin
Sevastopol National Technical University,
Sevastopol, Ukraine, 99053

**Abstract**

We present an asymptotic evaluation unitary formula for large argument values existing for defined class of functions. The asymptotic evaluation is obtained using only power series expansion coefficients of a function, what is a new result in power series analytic continuation theory.

# Introduction

In present the problem of obtaining information about function defined by power series depends significantly on its radius of convergence [1]. The universal approach to non-formal series is to compose and tabulate analytic continuations of it. The disadvantages of this way are large amount of intermediate coefficients and necessity of computing scheme stability analysis. If the series converges to an entire function in some cases we can research the sum with asymptotic methods [2]. The other most used method, suitable for formal series too, is to analyze asymptotic behaviour of a sequence of functional approximations composed for the function. This scheme differs much for different functions.

Some of used approximations are: method of Padé approximations [3], method of $G^3J$ - approximations [4,5], method of quadratic approximations [6], method of Levin approximations [7].

General results concerning modern state of Padé approximants application to the problem are described in paper [8]. Some other rational approximations of power series and ways of their application are described in [9,10].

In accordance with theory of analytic continuation of a complex variable function [11] all properties of a function analytic in some point are defined by its power expansion in this point. The central practical problem of the theory is research of function properties immediately on series coefficients given. Nonetheless the problem of research separate function behaviour is still actual. Some of the results are represented in papers [12 - 15]. In present paper let us formulate and prove the theorem concerning asymptotic evaluation of some class of functions with unitary formula containing power series coefficients of these functions.

The main approach of present paper is to prove that formulas for two first coefficients of asymptotic expansion in negative power series obtained in [16] under conditions connected with coefficients of Taylor expansion in some point of variable can be applied to any function fulfilling other conditions connected with properties of the function in infinite increasing of variable.

## Notations

$z$ – complex variable, $z = x + i \cdot y$, $i = \sqrt{-1}$.
$D$ – complex plane region, $D = A \cup B$.
$A$ – interior of complex plane circle $(x-a)^2 + y^2 = b^2$;
$a$ – real value, $a \in (-\infty, \infty)$;
$b$ – real value, $b \in (0, \infty)$;
$B$ – complex plane region, $x \in [a, \infty)$, $y \in [-b, b]$;
$H$ – interior of complex plane circle $x^2 + y^2 = h^2$;
$h$ – real value, $h \in (2, \infty)$;
$x_0$ – real value, $x_0 \in [a, \infty)$.

## Asymptotic evaluation of some class of functions

Theorem 1. *Let $f(z)$ be a complex variable function, analytic in all points of region D. Then if there exists a complex variable function $v(z)$ analytic in all points of region H and coinciding for real positive values of variable with $f\left(\dfrac{1}{x} + x_0 - 1\right)$, then the asymptotic evaluation with finite values of coefficients exists*

$$f(x) = q'_0 + q'_1 \cdot \left(\frac{1}{x}\right) + O\left(\frac{1}{x^2}\right), \quad x \to \infty, \tag{1}$$

*where*

$$q'_0 = \lim_{m \to \infty} \sum_{n=0}^{m} \binom{m}{n} \cdot c_n, \tag{2}$$

$$q'_1 = \lim_{m \to \infty} \sum_{n=1}^{m} \left\{ \binom{m}{n+1} - m \cdot \binom{m}{n} \right\} \cdot c_n,$$

*where $c_n$ are coefficients of Taylor expansion of $f(x)$ with center of convergence $(x = x_0)$.*

Proof.

Lemma 1.1. *Coefficients $q'_0$ and $q'_1$ of formal expansion $f(x)$ in series*

$$f(x) = \sum_{n=0}^{\infty} \frac{q'_n}{(x - x_0 + 1)^n} \qquad (3)$$

*are invariant with respect to $x_0$.*

Lemma 1.2. *The asymptotic evaluation exists*

$$f(x) = q'_0 + \frac{q'_1}{x} + O\left(\frac{1}{x^2}\right), \quad x \to \infty. \qquad (4)$$

Lemma 1.3. *In approximation of $f(x)$ by sequence*

$$R(x) = q'_{0,m} + \sum_{n=1}^{m} \frac{q'_{n,m}}{(x - x_0 + 1)^n}, \qquad (5)$$

*composed under condition of equality of first $(m+1)$ coefficients of power series expansions of $f(x)$ and $R(x)$ with center of convergence $x = x_0$, the coefficient of zero term is defined by formulas*

$$q'_{0,m} = \sum_{n=0}^{m} \binom{m}{n} \cdot c_n, \qquad (6)$$

$$q'_{1,m} = \sum_{n=1}^{m} \left\{ \binom{m}{n+1} - m \cdot \binom{m}{n} \right\} \cdot c_n. \qquad (7)$$

Lemma 1.4. *In unlimited increasing of approximation (5) dimension m coefficients of zero term of expansion (3) and approximation (5) coincide.*

Lemma 1.5. *In unlimited increasing of approximation (5) dimension m coefficients of first term of expansion (2) and approximation (5) coincide.*

## Proofs of the Lemmas

Lemma 1.1

Sufficient condition of existence of expansion $f(x)$ in series (3) is existence of $f(x)$ and all its derivatives for $\left(\frac{1}{x - x_0 + 1} = 0\right)$.

Let us note

$$\bar{x} = \frac{1}{x - x_0 + 1},$$

then

$$x = \frac{1 + (x_0 - 1) \cdot \bar{x}}{\bar{x}}.$$

Coefficients of expansion (3) are defined by Taylor formula

$$q'_n = \frac{1}{n!} \cdot \frac{\partial^n f\left(\frac{1+(x_0-1)\cdot \bar{x}}{\bar{x}}\right)}{\partial \bar{x}^n}\Bigg|_{\bar{x}=0}. \tag{8}$$

Let us consider the expansion of $f(x)$ in series

$$f(x) = \sum_{n=0}^{\infty} \frac{q''_n}{(x-x_1+1)^n}, \tag{9}$$

where $x_1 > x_0$.

Let us note

$$\bar{\bar{x}} = \frac{1}{x-x_1+1},$$

then

$$x = \frac{1+(x_1-1)\cdot \bar{\bar{x}}}{\bar{\bar{x}}},$$

$$\bar{\bar{x}} = \frac{\bar{x}}{1+(x_0-x_1)\cdot \bar{x}}, \tag{10}$$

$$\bar{x} = \frac{\bar{\bar{x}}}{1+(x_1-x_0)\cdot \bar{\bar{x}}}. \tag{11}$$

We have

$$q''_n = \frac{1}{n!} \cdot \frac{\partial^n f\left(\frac{1+(x_1-1)\cdot \bar{\bar{x}}}{\bar{\bar{x}}}\right)}{\partial \bar{\bar{x}}^n}\Bigg|_{\bar{\bar{x}}=0}. \tag{12}$$

Let us consider coefficients $q''_n$ for different $n$.

1. $n = 0$.

$$q''_0 = \lim_{\bar{\bar{x}} \to 0} f\left(\frac{1+(x_1-1)\cdot \bar{\bar{x}}}{\bar{\bar{x}}}\right). \tag{13}$$

We have

$$\lim_{\bar{x}\to 0} \frac{\bar{\bar{x}}}{\bar{x}} = \lim_{\bar{x}\to 0} \frac{1}{1+(x_0-x_1)\cdot \bar{x}} \equiv 1,$$

then

$$q''_0 = \lim_{\bar{x}\to 0} f\left(\frac{1+(x_0-1)\cdot \bar{x}}{\bar{x}}\right). \tag{14}$$

In accordance with (8),

$$q'_0 = \lim_{\bar{x}\to 0} f\left(\frac{1+(x_0-1)\cdot \bar{x}}{\bar{x}}\right). \tag{15}$$

As right parts of (14) and (15) coincide for any $x_1$, $q'_0$ is invariant in respect with $x_0$.

2. $n=1$.

$$q''_1 = \lim_{\bar{\bar{x}} \to 0} \frac{\partial f\left(\frac{1+(x_1-1)\cdot \bar{\bar{x}}}{\bar{\bar{x}}}\right)}{\partial \bar{\bar{x}}}. \tag{16}$$

In accordance with property of derivation operation we have

$$\frac{\partial f\left(\frac{1+(x_1-1)\cdot \bar{\bar{x}}}{\bar{\bar{x}}}\right)}{\partial \bar{\bar{x}}} = \frac{\partial f\left(\frac{1+(x_1-1)\cdot \bar{\bar{x}}}{\bar{\bar{x}}}\right)}{\partial \bar{x} \cdot \frac{\partial \bar{\bar{x}}}{\partial \bar{x}}},$$

or

$$\frac{\partial f\left(\frac{1+(x_1-1)\cdot \bar{\bar{x}}}{\bar{\bar{x}}}\right)}{\partial \bar{\bar{x}}} = \frac{\partial f\left(\frac{1+(x_1-1)\cdot \bar{\bar{x}}}{\bar{\bar{x}}}\right)}{\partial \bar{x}} \cdot \frac{\partial \bar{x}}{\partial \bar{\bar{x}}}.$$

Let us consider derivatives $\frac{\partial^n \bar{x}}{\partial \bar{\bar{x}}^n}$ for $\bar{\bar{x}} \to 0$.

In accordance with (11),

$$\bar{x} = \bar{\bar{x}} \cdot \sum_{n=0}^{\infty} (-1)^n \cdot (x_1 - x_0)^n \cdot \bar{\bar{x}}^n,$$

or

$$\bar{x} = \sum_{n=0}^{\infty} (-1)^n \cdot (x_1 - x_0)^n \cdot \bar{\bar{x}}^{n+1}. \tag{17}$$

(17) can be written as

$$\bar{x} = \sum_{n=0}^{\infty} d_n \cdot \bar{\bar{x}}^n,$$

where

$$d_n = \frac{1}{n!} \cdot \frac{\partial^n \bar{x}}{\partial \bar{\bar{x}}^n}\bigg|_{\bar{\bar{x}}=0}.$$

Then we have

$$\frac{\partial^n \bar{x}}{\partial \bar{\bar{x}}^n}\bigg|_{\bar{\bar{x}}=0} = n! \cdot d_n. \tag{18}$$

In accordance with (18),

$$\lim_{\bar{\bar{x}} \to 0} \frac{\partial^n \bar{x}}{\partial \bar{\bar{x}}^n} = n! \cdot (-1)^{n-1} \cdot (x_1 - x_0)^{n-1}. \tag{19}$$

In particulars,

$$\lim_{\bar{\bar{x}} \to 0} \frac{\partial \bar{x}}{\partial \bar{\bar{x}}} = 1! \cdot (-1)^0 \cdot (x_1 - x_0)^0 \equiv 1. \tag{20}$$

In accordance with (20),

$$\lim_{\bar{\bar{x}} \to 0} \frac{\partial f\left(\frac{1+(x_1-1)\cdot \bar{\bar{x}}}{\bar{\bar{x}}}\right)}{\partial \bar{x}} \cdot \frac{\partial \bar{x}}{\partial \bar{\bar{x}}} = \lim_{\bar{\bar{x}} \to 0} \frac{\partial f\left(\frac{1+(x_1-1)\cdot \bar{\bar{x}}}{\bar{\bar{x}}}\right)}{\partial \bar{x}}. \qquad (21)$$

We have

$$\lim_{\bar{\bar{x}} \to 0} \frac{\bar{x}}{\bar{\bar{x}}} = 1,$$

so (16) can be transformed to

$$q_1'' = \lim_{\bar{x} \to 0} \frac{\partial f\left(\frac{1+(x_0-1)\cdot \bar{x}}{\bar{x}}\right)}{\partial \bar{x}}. \qquad (22)$$

Right part (21) can be immediately obtained from (8) for $(n=1)$ and any $x_1$. Then $q_1'$ is invariant in respect with $x_0$.

We notify that coefficients of higher terms are not invariant in respect with $x_0$.

Lemma 1.2.

Under Theorem 1 hypothesis $v(z)$ is analytic in region $H$ and coincides for real positive values of variable with $f\left(\frac{1}{x}+x_0-1\right)$. Then the expansion exists

$$f\left(\frac{1}{\bar{x}}+x_0-1\right) = \sum_{n=0}^{\infty} q_n' \cdot \bar{x}^n, \qquad (23)$$

where

$$\bar{x} = \frac{1}{x-x_0+1}$$

and the series

$$S = \sum_{n=0}^{\infty} q_n' \cdot \bar{x}^n \qquad (24)$$

has radius of convergence $R' = h$.

Then the expansion of $f(x)$ in series (3) converges while

$$\bar{x} < h$$

or

$$x > x_0 - 1 + \frac{1}{h}.$$

For converging power series there exists [2] the asymptotic evaluation leading to

$$v(x) = \sum_{n=0}^{m} q_n' \cdot x^n + O(x^{m+1}), \quad x \in (-h, h), \qquad (25)$$

and

$$f(x) = q'_0 + \frac{q'_1}{x - x_0 + 1} + O\left(\frac{1}{(x - x_0 + 1)^2}\right), \quad x > x_0 - 1 + \frac{1}{h}. \tag{26}$$

In accordance with proved invariant values of $q'_0$ and $q'_1$ in respect with $x_0$ and the limit

$$\lim_{x \to \infty} \frac{x}{x + a} = 1, \tag{27}$$

we have

$$f(x) = q_0 + \frac{q_1}{x} + O\left(\frac{1}{x^2}\right), \quad x \to \infty. \tag{28}$$

Lemma 1.3.

Let us expand $R(x)$ (5) in powers of $(x - x_0)$:

$$R(x) = q'_{0,m} + \sum_{k=1}^{m} q'_{k,m} \cdot \sum_{n=0}^{\infty} (-1)^n \cdot \binom{k+n-1}{n} \cdot (x - x_0)^n. \tag{29}$$

After equating the coefficients for equal powers of variable we have

$$\begin{cases} c_0 = \sum_{k=0}^{m} q'_{k,m} \\ c_n = (-1)^n \cdot \sum_{k=1}^{m} q'_{k,m} \cdot \binom{k+n-1}{n} \end{cases} \tag{30}$$

In accordance with (5) and (30) the system of linear algebraic equations connecting coefficients $c$ and $q'_{,m}$ is

$$\begin{cases} c_0 = \sum_{k=0}^{m} q'_{k,m} \\ \cdots \\ c_n = (-1)^n \cdot \sum_{k=1}^{m} q'_{k,m} \cdot \binom{k+n-1}{n} \\ \cdots \\ c_m = (-1)^m \cdot \sum_{k=1}^{m} q'_{k,m} \cdot \binom{k+m-1}{m} \end{cases} \tag{31}$$

Let us transform the system (31) to a triangle form. For matrix defined by right part of (31) Gauss's method can be applied next way.

Definition:

*Transformation of k-th order of countable set R transforms it to countable set R' so that*

$$\begin{cases} R'_i = R_i & \text{for } i = 0,...,k \\ R'_i = R_i + R_{i-1} & \text{for } i > k \end{cases} \tag{32}$$

Proposition 1.
*If R is a countable set with elements defined by equations*

$$\begin{cases} R_i = c_0 & \text{for } i = 0 \\ R_i = \sum_{s=0}^{i-1} \binom{i-1}{s} \cdot c_{i-s} & \text{for } i = 1,...,k, \\ R_i = \sum_{s=0}^{k-1} \binom{k-1}{s} \cdot c_{i-s} & \text{for } i > k \end{cases} \tag{33}$$

*then transformation of k-th order transforms it to countable set $R'$ with elements*

$$\begin{cases} R'_i = c_0 & \text{for } i = 0 \\ R'_i = \sum_{s=0}^{i-1} \binom{i-1}{s} \cdot c_{i-s} & \text{for } i = 1,...,k+1 \\ R'_i = \sum_{s=0}^{k} \binom{k}{s} \cdot c_{i-s} & \text{for } i > k+1 \end{cases} \tag{34}$$

Proof.
Obviously for $i = 0,...,k$ Proposition 1 is correct.

Let us consider $i = k + 1$. In accordance with (32), element of $R'$ should be equal the sum of $R_i$ and $R_{i-1}$.

We have next equations. From (34):

$$R'_i = R'_{k+1} = \sum_{s=0}^{k} \binom{k}{s} \cdot c_{k+1-s}. \tag{35}$$

From (33):

$$R_i = R_{k+1} = \sum_{s=0}^{k-1} \binom{k-1}{s} \cdot c_{k+1-s}. \tag{36}$$

From (33):

$$R_{i-1} = R_k = \sum_{s=0}^{k-1} \binom{k-1}{s} \cdot c_{k-s}. \tag{37}$$

Then

$$R_i + R_{i-1} = \sum_{s=0}^{k-1} \binom{k-1}{s} \cdot c_{k+1-s} + \sum_{s=0}^{k-1} \binom{k-1}{s} \cdot c_{k-s}.$$

Let us note $k - s = k + 1 - s'$, then $s = s' - 1$ or $s' = s + 1$, for $s = 0$ $s' = 1$, for $s = k - 1$ $s' = k$, then we have

$$R_i + R_{i-1} = \sum_{s=0}^{k-1}\binom{k-1}{s}\cdot c_{k+1-s} + \sum_{s'=1}^{k-1}\binom{k-1}{s'-1}\cdot c_{k+1-s'},$$

or, after replacement in second additive $s'$ by $s$,

$$R_i + R_{i-1} = \sum_{s=0}^{k-1}\binom{k-1}{s}\cdot c_{k+1-s} + \sum_{s=1}^{k-1}\binom{k-1}{s-1}\cdot c_{k+1-s}. \tag{38}$$

The identities exist

$$\binom{k-1}{0} \equiv 1 \text{ и } \binom{k-1}{k-1} \equiv 1,$$

so (38) can be transformed to

$$R_i + R_{i-1} = c_{k-1} + \sum_{s=1}^{k-1}\left\{\binom{k-1}{s}+\binom{k-1}{s-1}\right\}\cdot c_{k+1-s} + c_1. \tag{39}$$

As the identity takes place [17]

$$\binom{k-1}{s}+\binom{k-1}{s-1} \equiv \binom{k}{s},$$

we have

$$R_i + R_{i-1} = c_{k+1} + \sum_{s=1}^{k-1}\binom{k}{s}\cdot c_{k+1-s} + c_1 = \sum_{s=0}^{k}\binom{k}{s}\cdot c_{k+1-s}, \tag{40}$$

what coincides with (35).

For $i = k+1$ Proposition 1 is correct.

Let us consider $i > k+1$. Under Definition (32), the element of $R'_i$ should be equal the sum of $R_i$ and $R_{i-1}$.

Then we have next equations. From (34):

$$R'_i = \sum_{s=0}^{k}\binom{k}{s}\cdot c_{i-s}.$$

From (33):

$$R_i = \sum_{s=0}^{k-1}\binom{k-1}{s}\cdot c_{i-s}.$$

From (33):

$$R_{i-1} = \sum_{s=0}^{k-1}\binom{k-1}{s}\cdot c_{i-1-s}.$$

Then

$$R_i + R_{i+1} = \sum_{s=0}^{k-1}\binom{k-1}{s}\cdot c_{i-s} + \sum_{s=0}^{k-1}\binom{k-1}{s}\cdot c_{i-1-s}.$$

Let us note $i-1-s = i-s'$, or $i-(1+s) = i-s'$, then $s' = s+1$, $s = s'-1$, for $s = 0$ $s' = 1$, for $s = k-1$ $s' = k$, then we have

$$R_i + R_{i-1} = \sum_{s=0}^{k-1}\binom{k-1}{s}\cdot c_{i-s} + \sum_{s'=1}^{k}\binom{k-1}{s'-1}\cdot c_{i-s'},$$

or, after replacement $s'$ by $s$ in second additive of right part,

$$R_i + R_{i-1} = \sum_{s=0}^{k-1} \binom{k-1}{s} \cdot c_{i-s} + \sum_{s=1}^{k} \binom{k-1}{s'-1} \cdot c_{i-s}. \tag{41}$$

The identities exist
$$\binom{k-1}{0} \equiv \binom{k-1}{k-1} \equiv 1,$$
$$\binom{k-1}{s} + \binom{k-1}{s-1} \equiv \binom{k}{s}.$$

Then (41) can be transformed to
$$R_i + R_{i-1} = c_i + \sum_{s=1}^{k-1} \binom{k}{s} \cdot c_{i-s} + c_{i-k},$$

or
$$R_i + R_{i-1} = \sum_{s=0}^{k} \binom{k}{s} \cdot c_{i-s}, \tag{42}$$

what coincides with third equation of (24).

Then Proposition 1 is correct for $i > k+1$.

Proposition 1 is proved.

Proposition 2.

*If R is a countable set with elements defined by equations*
$$\begin{cases} R_i = (-1)^i \cdot \binom{j}{j-i} & \text{при } i=0,...,k \\ R_i = (-1)^i \cdot \binom{i+j-k}{j-k} & \text{при } i>k \end{cases}, \tag{43}$$

*then transformation of k-th order transforms it to countable set $R'$ with elements*
$$\begin{cases} R'_i = (-1)^i \cdot \binom{j}{j-i} & \text{при } i=0,...,k+1 \\ R'_i = (-1)^i \cdot \binom{i+j-k-1}{j-k-1} & \text{при } i>k+1 \end{cases}. \tag{44}$$

Proof.

Obviously for $i = 0,...,k$ Proposition 2 is correct.

Let us consider $i = k+1$. Under Definition (32) the element of $R'_i$ should be equal the sum of $R_i$ and $R_{i-1}$.

We have next equations. From (44):
$$R'_i = R'_{k+1} = (-1)^{k+1} \cdot \binom{j}{j-k-1}. \tag{45}$$

From (43):
$$R_i = R_{k+1} = (-1)^{k+1} \cdot \binom{j+1}{j-k}. \tag{46}$$

From (43):
$$R_{i-1} = R_k = (-1)^k \cdot \binom{j}{j-k}. \tag{47}$$

Then
$$R_i + R_{i-1} = (-1)^{k+1} \cdot \left\{ \binom{j+1}{j-k} - \binom{j}{j-k} \right\}.$$

In accordance with property of binomial coefficients [17],
$$\binom{a}{b} + \binom{a}{b-1} \equiv \binom{a+1}{b},$$

then
$$\binom{a+1}{b} - \binom{a}{b} = \binom{a}{b-1},$$

or
$$\binom{j+1}{j-k} - \binom{j}{j-k} = \binom{j}{j-k-1}.$$

Then
$$R_i + R_{i-1} = (-1)^{k+1} \cdot \binom{j}{j-k-1}, \tag{48}$$

what coincides with (45). For $i = k+1$ Proposition 2 is correct.

Let us consider $i > k+1$. Under Definition (32) the element of $R'_i$ should be equal the sum of $R_i$ and $R_{i-1}$.

In accordance with (43),
$$R_i = (-1)^i \cdot \binom{i+j-k}{j-k},$$
$$R_{i-1} = (-1)^{i-1} \cdot \binom{i+j-k-1}{j-k}.$$

In accordance with (44),
$$R'_i = (-1)^i \cdot \binom{i+j-k-1}{j-k-1}.$$

Then
$$R_i + R_{i-1} = (-1)^i \cdot \left\{ \binom{i+j-k}{j-k} - \binom{i+j-k-1}{j-k} \right\}.$$

In accordance with the identity
$$\binom{a}{b} - \binom{a-1}{b} \equiv \binom{a-1}{b-1}$$

we have
$$R_i + R_{i-1} = (-1)^i \cdot \binom{i+j-k-1}{j-k-1}, \tag{49}$$

what coincides with second formula of (44). Then Proposition 2 is correct for $i > k+1$.

Proposition 2 is proved.

Let us note transformation of *k*-th order of countable set $R$ as $T_k(R)$. Then in notations of (32) we have
$$R' = T_k(R). \tag{50}$$

Definition:

*Sequential transformation of m-th order of countable set $R$ transforms it to countable set $R^{(m)}$ by processing of recurrent formula*
$$R^{(k)} = T_k\left(R^{(k-1)}\right) \tag{51}$$
*for $k = 1,...,m$, where $R^{(0)} = R$.*

Proposition 3.

*Elements of countable set $R^{(m)}$ are connected with elements of countable set $R$ by equations*

$$\begin{cases} R^{(m)}{}_0 = R_0 \\ R^{(m)}{}_i = \sum_{s=0}^{i-1}\binom{i-1}{s} \cdot R_{i-s} \quad \text{for} \quad i = 1,...,m+1. \\ R^{(m)}{}_i = \sum_{s=0}^{m}\binom{m}{s} \cdot R_{i-s} \quad \text{for} \quad i > m+1 \end{cases} \tag{52}$$

Proof.

For $m = 1$ $R^{(1)} = T_1\left(R^{(0)}\right)$, or $R' = T_1(R)$. In accordance with (32) we have next values of elements of $R^{(1)}$:
$$\begin{cases} \text{for } i = 0,...,1 \quad R'_i = R_i \\ \text{for } i > 1 \quad R'_i = R_i + R_{i-1} \end{cases}. \tag{53}$$

From (52) we have next equations
$$\begin{cases} R^{(1)}{}_0 = R_0 \\ \text{for } i = 1: \quad R^{(1)}{}_1 = \sum_{s=0}^{0}\binom{0}{s} \cdot R_{1-s} = R_1 \\ \text{for } i = 2: \quad R^{(1)}{}_2 = \sum_{s=0}^{1}\binom{1}{s} \cdot R_{2-s} = R_2 + R_1 \\ \text{for } i > 2: \quad R^{(1)}{}_i = \sum_{s=0}^{1}\binom{1}{s} \cdot R_{i-s} = R_i + R_{i-1} \end{cases},$$

what coincides with (53). For $m=1$ Proposition 3 is correct.

Let us use induction hypothesis. Let Proposition 3 be correct for some $m$ by assumption. If Proposition 3 is correct for $(m+1)$ then elements of $R^{(m+1)}$ should be connected with elements of $R$ by equations

$$\begin{cases} R^{(m+1)}{}_0 = R_0 \\ \text{for } i=1,\ldots,m+2 \quad R^{(m+1)}{}_i = \sum_{s=0}^{i-1} \binom{i-1}{s} \cdot R_{i-s} \\ \text{for } i>m+2 \quad R^{(m+1)}{}_i = \sum_{s=0}^{m+1} \binom{m+1}{s} \cdot R_{i-s} \end{cases} \qquad (54)$$

Under Definition (51), $R^{(m+1)}$ is a result of sequential processing of operation $R^{(k)} = T_k(R^{(k-1)})$ for $k=1,\ldots,m+1$. But $R^{(m)}$ is the result of the same operations for $k=1,\ldots,m$.

Then
$$R^{(m+1)} = T_{m+1}(R^{(m)}). \qquad (55)$$

Let us define $k=m+1$. Then, in accordance with truth of Proposition 3 for $m$, after replacement $R^{(m)}{}_i$ by $R_i$ and $R_{i-s}$ by $c_{i-s}$ we have equations of Proposition 1. After replacement in (34) $R'_i$ by $R^{(m+1)}{}_i$ and $c_{i-s}$ by $R_{i-s}$ we have (54).

Proposition 3 is proved.

Proposition 4.
*If R is a countable set with elements*
$$R_i = (-1)^i \cdot \binom{i+j-1}{j-1}, \qquad (56)$$
*where $j \geq 0$, then elements of countable set $R^{(m)}$ are defined by equations*

$$\begin{cases} \text{for } i=0,\ldots,m+1 \quad R^{(m)}{}_i = (-1)^i \cdot \binom{j}{j-i} \\ \text{for } i>m+1 \quad R^{(m)}{}_i = (-1)^i \cdot \binom{i+j-m-1}{j-m-1} \end{cases} \qquad (57)$$

Proof.

For $m=1$ $R^{(m)} = T_1(R_0)$, or $R' = T_1(R)$. In accordance with (32) and (56), the values of elements of $R$ are

$$\begin{cases} \text{for } i=0,\ldots,1 \quad R'_i = (-1)^i \cdot \binom{i+j-1}{j-1} \\ \text{for } i>1 \quad R'_i = R_i + R_{i-1} \end{cases}$$

Formula (56) is a special case of (43) for $k=1$ because

$$R_1 = (-1)^1 \cdot \binom{1+j-1}{j-1} = (-1)^1 \cdot \binom{j}{j-1} = (-1)^i \cdot \binom{j}{j-i}\bigg|_{i=1}.$$

In accordance with Proposition 2, transformation of first order transforms countable set $R$ with elements (56) to countable set $R'$ with elements (57), where $m = 1$. So, for $(m=1)$ Proposition 4 is correct.

Let us use induction hypothesis. Let Proposition 4 be correct for some $m > 1$ by assumption. If Proposition 4 is correct for $(m+1)$, then elements of $R^{(m+1)}$ should be defined by equations

$$\begin{cases} \text{for } i = 0,...,m+2 \quad R^{(m+1)}{}_i = (-1)^i \cdot \binom{j}{j-i} \\ \text{for } i > m+2 \quad R^{(m+1)}{}_i = (-1)^i \cdot \binom{i+j-m-2}{j-m-2} \end{cases} \quad (58)$$

In accordance with $R^{(m+1)} = T_{m+1}(R^{(m)})$ and truth of (57) for $m$ by assumption, after replacement in (57) $(m+1)$ by $k$ and $R^{(m)}{}_i$ by $R_i$, we have (43). In accordance with Proposition 2, after $k$-th order transformation of countable set (43) we will obtain countable set (44). After replacement in (44) $k$ by $(m+1)$ and $R'_i$ by $R^{(m+1)}{}_i$ we have (58).

Proposition 4 is proved.

Proposition 5.
*If operation of addition of countable sets is defined as*
$$U + V = R, \quad (59)$$
*where*
$$R_i = U_i + V_i, \quad (60)$$
*then*
$$(U+V)^{(m)} = U^{(m)} + V^{(m)}. \quad (61)$$

Proof.
Let us note $U + V = R$. Let us consider $R^{(m)}$. From (52) we have

$$\begin{cases} R^{(m)}{}_0 = R_0 = U_0 + V_0 \\ R^{(m)}{}_i = \sum_{s=0}^{i-1} \binom{i-1}{s} \cdot R_{i-s} = \sum_{s=0}^{i-1} \binom{i-1}{s} \cdot (U_{i-s} + V_{i-s}) \text{ for } i = 1,...,m+1, \\ R^{(m)}{}_i = \sum_{s=0}^{m} \binom{m}{s} \cdot R_{i-s} = \sum_{s=0}^{m} \binom{m}{s} \cdot (U_{i-s} + V_{i-s}) \text{ for } i > m+1 \end{cases} \quad (63)$$

$$\begin{cases} U^{(m)}{}_0 = U_0 \\ U^{(m)}{}_i = \sum_{s=0}^{i-1} \binom{i-1}{s} \cdot U_{i-s} \quad \text{for} \quad i=1,\ldots,m+1, \\ U^{(m)}{}_i = \sum_{s=0}^{m} \binom{m}{s} \cdot U_{i-s} \quad \text{for} \quad i > m+1 \end{cases}$$

$$\begin{cases} V^{(m)}{}_0 = V_0 \\ V^{(m)}{}_i = \sum_{s=0}^{i-1} \binom{i-1}{s} \cdot V_{i-s} \quad \text{for} \quad i=1,\ldots,m+1. \\ V^{(m)}{}_i = \sum_{s=0}^{m} \binom{m}{s} \cdot V_{i-s} \quad \text{for} \quad i > m+1 \end{cases}$$

Then, in accordance with (60), elements of $\left(U^{(m)} + V^{(m)}\right)$ are defined by equations

$$\begin{cases} \left(U^{(m)} + V^{(m)}\right)_0 = U_0 + V_0 \\ \left(U^{(m)} + V^{(m)}\right)_i = \sum_{s=0}^{i-1} \binom{i-1}{s} \cdot (U_{i-s} + V_{i-s}) \quad \text{for} \quad i=1,\ldots,m+1, \\ \left(U^{(m)} + V^{(m)}\right)_i = \sum_{s=0}^{m} \binom{m}{s} \cdot (U_{i-s} + V_{i-s}) \quad \text{for} \quad i > m+1 \end{cases}$$

what coincides with (62).

Proposition 5 is proved.

Let us write system (31) next way:
$$C = Q'_0 + \sum_{j=1}^{m} Q'_j, \tag{64}$$

where $C$ is a countable set with elements
$$\begin{cases} C_i = c_i \quad \text{for} \quad i \leq m \\ C_i = 0 \quad \text{for} \quad i > m \end{cases},$$

$Q'_0$ is a countable set with elements
$$\begin{cases} Q'_{00} = q'_{0,m} \\ Q'_{0i} = 0 \end{cases},$$

$Q'_j$ is a countable set with elements

$$\begin{cases} Q'_{ji} = (-1)^i \cdot q'_{j,m} \cdot \binom{i+j-1}{j-1} & \text{for } i \leq m \\ Q'_{ji} = 0 & \text{for } i > m \end{cases}.$$

Let us apply sequential transformation of $m$-th order to left and right parts of (64). After transformation we have:

$$C^{(m)} = Q_0'^{(m)} + \sum_{j=1}^{m} Q_j'^{(m)}. \tag{65}$$

In accordance with (52) and (56) we shall obtain next equations for nonzero elements of $C^{(m)}$:

$$\begin{cases} C^{(m)}{}_0 = \sum_{j=0}^{m} q'_{j,m} \\ \dots\dots\dots\dots\dots\dots\dots\dots\dots\dots \\ C^{(m)}{}_i = (-1)^i \cdot \sum_{j=1}^{m} \binom{j}{j-i} \cdot q'_{j,m} \\ \dots\dots\dots\dots\dots\dots\dots\dots\dots\dots \\ C^{(m)}{}_m = (-1)^m \cdot \sum_{j=1}^{m} \binom{j}{j-m} \cdot q'_{j,m} \end{cases} \tag{66}$$

In accordance with identities

$$\binom{a}{a-b} \equiv \binom{a}{b}$$

and

$$\binom{a}{b} \equiv 0$$

for $b < 0$ we have finally:

$$\begin{cases} C^{(m)}{}_0 = \sum_{j=0}^{m} q'_{j,m} \\ C^{(m)}{}_1 = -\sum_{j=1}^{m} j \cdot q'_{j,m} \\ \dots\dots\dots\dots\dots\dots\dots\dots\dots\dots \\ C^{(m)}{}_i = (-1)^i \cdot \sum_{j=i}^{m} \binom{j}{i} \cdot q'_{j,m} \\ \dots\dots\dots\dots\dots\dots\dots\dots\dots\dots \\ C^{(m)}{}_{m-1} = (-1)^{m-1} \cdot \{q'_{m-1,m} + m \cdot q'_{m,m}\} \\ C^{(m)}{}_m = (-1)^m \cdot q'_{m,m} \end{cases} \tag{67}$$

From linear algebra [18] is known that determinant of a triangle matrix is equal the product of elements of main diagonal. Countable set with finite quantity of nonzero elements can be considered as a vector of dimension $m$ so (67) can be written as

$$\overline{C}^{(m)} = (A) \cdot \overline{q}_m, \qquad (68)$$

where

$$\begin{cases} \overline{C}^{(m)}{}_0 = c_0 \\ \overline{C}^{(m)}{}_1 = c_1 \\ \overline{C}^{(m)}{}_i = \sum_{s=0}^{i-1} \binom{i-1}{s} \cdot c_{i-s} \end{cases},$$

$(A)$ is a square matrix of dimension $m$ with elements

$$a_{i,j} = (-1)^i \cdot \binom{j}{i}, \; i = 0,\ldots,m, \; j = 0,\ldots,m,$$

or

$$(A) = \begin{pmatrix} 1 & 1 & 1 & 1 & 1 & . & 1 \\ 0 & -1 & -2 & -3 & -4 & . & -m \\ 0 & 0 & 1 & 3 & 6 & . & \binom{m}{2} \\ 0 & 0 & 0 & -1 & -4 & . & -\binom{m}{3} \\ 0 & 0 & 0 & 0 & 1 & . & \binom{m}{4} \\ . & . & . & . & . & . & . \\ 0 & 0 & 0 & 0 & 0 & . & 1 \end{pmatrix}. \qquad (69)$$

Then determinant of (69) is equal

$$\mathrm{Det}(A) = \prod_{i=0}^{m}(-1)^i. \qquad (70)$$

(70) is nonzero for any nonnegative $m$, then in accordance with Cronecker-Capelli theorem, system (67) has unique solution

$$\overline{q}_m = (A)^{-1} \cdot \overline{C}^{(m)}. \qquad (71)$$

Let us introduce the functions

$$U(x) = \sum_{n=0}^{m} \overline{C}^{(m)}{}_n \cdot x^n, \qquad (72)$$

$$V(x) = \sum_{n=0}^{m} q'_{n,m} \cdot x^n. \qquad (73)$$

The identity exists

$$\frac{1}{k!} \cdot \frac{d^k}{dx^k} \cdot \sum_{n=0}^{m} a_n \cdot x^n \equiv \sum_{n=k}^{m} \binom{n}{k} \cdot x^{n-k}. \tag{74}$$

In accordance with (67) and (74) we have

$$C^{(m)}{}_i = (-1)^i \cdot \frac{1}{i!} \cdot \left\{ \frac{d^i}{dx^i} \cdot \sum_{j=0}^{m} q'_{j,m} \cdot x^j \right\} \bigg|_{x=1}, \tag{75}$$

or, in accordance with (73),

$$C^{(m)}{}_i = (-1)^i \cdot \frac{1}{i!} \cdot \frac{\partial^i V(x)}{\partial x^i} \bigg|_{x=1}. \tag{76}$$

As $V(x)$ is a polynomial, then there exists an expansion of $V(x)$ in powers of $(x-1)$:

$$V(x) = \sum_{n=0}^{m} v_n \cdot (x-1)^n, \tag{77}$$

where, in accordance with Taylor's formula,

$$v_n = \frac{1}{n!} \cdot \frac{\partial^n v(x)}{\partial x^n} \bigg|_{x=1}.$$

Then, in accordance with (76),

$$v_n = (-1)^n \cdot C^{(m)}{}_n. \tag{78}$$

In accordance with (78), we have

$$V(x) = \sum_{n=0}^{m} (-1)^n \cdot C^{(m)}{}_n \cdot (x-1)^n. \tag{79}$$

Using (79) and Taylor's formula, we have

$$q'_{n,m} = \frac{1}{n!} \cdot \frac{\partial^n V(x)}{\partial x^n} \bigg|_{x=0}. \tag{80}$$

In accordance with (79) and (74),

$$\frac{1}{n!} \cdot \frac{\partial^n V(x)}{\partial x^n} = \sum_{k=n}^{m} (-1)^k \cdot \binom{k}{n} \cdot (x-1)^{k-n} \cdot C^{(m)}{}_n. \tag{81}$$

For $(x=0)$

$$V(0) = \sum_{n=0}^{m} C^{(m)}{}_n, \tag{82}$$

$$\frac{1}{k!} \cdot \frac{\partial^n V(x)}{\partial x^n} \bigg|_{x=0} = \sum_{n=k}^{m} (-1)^{2n-k} \cdot \binom{n}{k} \cdot C^{(m)}{}_n = (-1)^k \cdot \sum_{n=k}^{m} \binom{n}{k} \cdot C^{(m)}{}_n. \tag{83}$$

In accordance with (80), (82), (83), the required formulas are

$$\begin{cases} q'_{0,m} = \sum_{n=0}^{m} C^{(m)}{}_n \\ q'_{1,m} = -\sum_{n=1}^{m} n \cdot C^{(m)}{}_n \\ q'_{k,m} = (-1)^k \cdot \sum_{n=k}^{m} \binom{n}{k} \cdot C^{(m)}{}_n \end{cases} \quad , \tag{84}$$

where

$$\begin{cases} C^{(m)}{}_0 = c_0 \\ C^{(m)}{}_n = \sum_{s=0}^{n-1} \binom{n-1}{s} \cdot c_{n-s} \end{cases} . \tag{85}$$

(84) can be transformed to

$$\overline{q}'_m = (A) \cdot \overline{C}^{(m)}. \tag{86}$$

But, in accordance with (71),

$$\overline{q}'_m = (A)^{-1} \cdot \overline{C}^{(m)}.$$

So for matrix (69) next equation is correct

$$(A)^{-1} = (A). \tag{87}$$

Right parts of (85) and (33)[16] coincides, so, using the deduction of formulas (126) [16], (84) can be transformed to

$$\begin{cases} q'_{0,m} = \sum_{s=0}^{m} c_s \cdot \binom{m}{s} \\ q'_{1,m} = -\sum_{s=1}^{m} c_s \cdot \left\{ m \cdot \binom{m}{s} - \binom{m}{s+1} \right\} \\ q'_{k,m} = (-1)^k \cdot \sum_{s=1}^{m} c_s \cdot \sum_{n=0}^{k} (-1)^n \cdot \binom{m-n}{k-n} \cdot \binom{m}{s+n} \end{cases} . \tag{88}$$

After replacement in first two formulas *s* by *n* we have (6) and (7).

Lemma 1.4.

Let us replace in system (31) $q'_{k,m}$ by $q'_k$, the difference of each equation after replacement let us note as $R_{n,m}$, then we have

$$\begin{cases} c_0 = \sum_{k=0}^{m} q'_k + R_{0,m} \\ \cdots\cdots\cdots\cdots\cdots\cdots\cdots\cdots\cdots\cdots\cdots\cdots \\ c_n = (-1)^n \cdot \sum_{k=1}^{m} q'_k \cdot \binom{k+n-1}{n} + R_{n,m} \\ \cdots\cdots\cdots\cdots\cdots\cdots\cdots\cdots\cdots\cdots\cdots\cdots \\ c_m = (-1)^m \cdot \sum_{k=1}^{m} q'_k \cdot \binom{k+m-1}{m} + R_{m,m} \end{cases} \quad (89)$$

In accordance with (88) and (89) we have

$$q'_0 = \sum_{n=0}^{m} \binom{m}{n} \cdot (c_n - R_{n,m}). \quad (90)$$

In infinite increasing of dimension $m$ of approximation (5) we have

$$q'_0 = \lim_{m\to\infty} \sum_{n=0}^{m} \binom{m}{n} \cdot (c_n - R_{n,m}),$$

or

$$q'_0 = \lim_{m\to\infty} \sum_{n=0}^{m} \binom{m}{n} \cdot c_n - \lim_{m\to\infty} \sum_{n=0}^{m} \binom{m}{n} \cdot R_{n,m}. \quad (91)$$

In accordance with (88) let us transform (91) to

$$\lim_{m\to\infty} q'_{0,m} = q'_0 + \lim_{m\to\infty} \sum_{n=0}^{m} \binom{m}{n} \cdot R_{n,m}. \quad (92)$$

Then the proof of Lemma 1.4 can be reduced to the proof of proposition

$$\lim_{m\to\infty} \sum_{n=0}^{m} \binom{m}{n} \cdot R_{n,m} = 0 \quad (93)$$

under Theorem 1 hypothesis.

Let us consider $R_{n,m}$.

In accordance with (1), (3) and (89) we have

$$\begin{cases} R_{0,m} = \sum_{k=m+1}^{\infty} q'_k \\ R_{n,m} = (-1)^n \cdot \sum_{k=m+1}^{\infty} q'_k \cdot \binom{k+n-1}{n} \end{cases} \quad (94)$$

Let us note

$$S = \sum_{n=0}^{m} \binom{m}{n} \cdot R_{n,m}.$$

In accordance with (94) we have

$$S = \sum_{n=0}^{m} \binom{m}{n} \cdot (-1)^n \cdot \sum_{k=m+1}^{\infty} q'_k \cdot \binom{k+n-1}{n}. \quad (95)$$

Let us consider the convergence of series of (94) right part. Under Theorem 1 conditions, $v(z)$ is a function analytic in region $H$, then radius of convergence of series (24) $R > 2$. As a power series converges absolutely inside its disk of convergence, then for $(x = 1)$ series (24) converges absolutely. Finite quantity of terms do not affect on the convergence of series, so right part of first equation of (94) also converges absolutely.

To research the convergence of other series of (94) let us consider the series obtained by $(m+1)$-th derivation of (24). In accordance with Abel's theorem [2] this series also converges absolutely inside disk of convergence of (24), so for $(x = 1)$.

Then we have

$$\frac{d^{m+1}}{dx^{m+1}} \cdot \sum_{k=0}^{\infty} q'_k \cdot x^k = (m+1)! \sum_{k=m+1}^{\infty} \binom{k}{m+1} \cdot q'_k \cdot x^{k-m-1}. \qquad (96)$$

Proposition 6.

*For nonnegative integers $k, n, m$, $k > m+1$, $n \le m$ we have*

$$(m+1)! \binom{k}{m+1} > \binom{k+n-1}{n}. \qquad (97)$$

Proof.

Let us define

$$S_1 = (m+1)! \binom{k}{m+1}, \quad S_2 = \binom{k+n-1}{n}.$$

In accordance with properties of binomial coefficients [17] we have

$$S_1 = \frac{k!}{(k-m-1)!}, \qquad (98)$$

$$S_2 = \frac{(k+n-1)!}{n! (k-1)!}. \qquad (99)$$

Then Proposition 6 is equal the proposition

$$k! \cdot n! \prod_{i=1}^{m} (k-i) > (k+n-1)!,$$

or

$$n! \prod_{i=1}^{m} (k-i) > \prod_{i=1}^{n-1} (k+i) \qquad (100)$$

under hypothesis of Proposition 6.

For $(n \le m)$ (100) can be transformed to

$$\left( \prod_{i=1}^{n-1} (i+1) \right) \cdot \left( \prod_{i=1}^{n-1} (k-i) \right) \cdot \prod_{i=n}^{m} (k-i) > \prod_{i=1}^{n-1} (k+i),$$

or

$$\left( \prod_{i=1}^{n-1} (i+1) \cdot (k-i) \right) \cdot \prod_{i=n}^{m} (k-i) > \prod_{i=1}^{n-1} (k+i). \qquad (101)$$

Sufficient condition of correctness of (101) is correctness of system

$$\begin{cases} (i+1) \cdot (k-i) > (k+i) & \text{for } i=1,\ldots,n-1 \\ \prod_{i=n}^{m}(k-i) \geq 1 \end{cases} \qquad (102)$$

under hypothesis of Proposition 6.

From first inequality of (102) we have

$$(k-2) \cdot i - i^2 > 0$$

then

$$i \in (0, k-2). \qquad (103)$$

As, in accordance with Proposition 6 hypothesis, $k > m+1$, then $i < m-1$. Then, for $n \leq m$ the fist inequality of (102) is correct within the entire interval $[1, n-1)$.

Under hypothesis of Proposition 6 the second inequality of (102) is correct automatically.

Proposition 6 is proved.

In accordance with (97) each term of series of (94) right part is less by absolute value than the term of the same order of series (96) for $(x=1)$, that converges absolutely. So all right parts of (94) converge absolutely.

Then we can change the order of summation of (94). Let us transform (94) to

$$S = \sum_{k=m+1}^{\infty} q'_k \cdot \sum_{n=0}^{m} (-1)^n \cdot \binom{m}{n} \cdot \binom{k+n-1}{n}. \qquad (104)$$

Let us consider the sum

$$S_1 = \sum_{n=0}^{m} (-1)^n \cdot \binom{m}{n} \cdot \binom{k+n-1}{n}. \qquad (105)$$

Proposition 7.

*For integers m, k, $m \geq 1$, $k \geq 0$, $m \geq k+1$ we have*

$$\sum_{n=0}^{k} (-1)^n \cdot \binom{m}{n} = (-1)^k \cdot \binom{m-1}{k} \qquad (106)$$

Proof.

For $(k=0)$ in accordance with (106) we have

Left part:

$$\sum_{n=0}^{0} (-1)^n \cdot \binom{m}{n} = \binom{m}{0} \equiv 1.$$

Right part:

$$(-1)^0 \cdot \binom{m-1}{0} \equiv 1.$$

Proposition 7 is correct.

For $(k=1)$ in accordance with (106) we have

Left part:

$$\sum_{n=0}^{1}(-1)^n \cdot \binom{m}{n} = \binom{m}{0} - \binom{m}{1} = 1 - m.$$

Right part:

$$(-1)^1 \cdot \binom{m-1}{1} \equiv -(m-1) = 1 - m.$$

Proposition 7 is correct.

Let us use induction hypothesis. Let Proposition 7 be correct for some integer $k > 1$ by assumption. If Proposition 7 is correct for $(k+1)$ then next equation should be correct

$$\sum_{n=0}^{k+1}(-1)^n \cdot \binom{m}{n} = (-1)^{k+1} \cdot \binom{m-1}{k+1}. \tag{107}$$

(107) can be transformed to

$$\sum_{n=0}^{k+1}(-1)^n \cdot \binom{m}{n} = \sum_{n=0}^{k}(-1)^n \cdot \binom{m}{n} + (-1)^{k+1} \cdot \binom{m}{k+1}. \tag{108}$$

As Proposition 7 is correct for $k$ by assumption, in accordance with (108) we have

$$\sum_{n=0}^{k+1}(-1)^n \cdot \binom{m}{n} = (-1)^k \cdot \binom{m-1}{k} + (-1)^{k+1} \cdot \binom{m}{k+1},$$

or

$$\sum_{n=0}^{k+1}(-1)^n \cdot \binom{m}{n} = (-1)^{k+1} \cdot \left\{\binom{m}{k+1} - \binom{m-1}{k}\right\}. \tag{109}$$

In accordance with "Pascal's triangle" we have

$$\binom{a+1}{b+1} - \binom{a}{b} \equiv \binom{a}{b+1}. \tag{110}$$

In accordance with (109) and (110),

$$\sum_{n=0}^{k+1}(-1)^n \cdot \binom{m}{n} = (-1)^{k+1} \cdot \binom{m-1}{k+1},$$

what coincides with (107). Then Proposition 7 is correct for $(k+1)$.

Proposition 7 is proved.

Let us introduce summation index $r = n+1$, then $n = r-1$, for $n=0$ $r=1$, for $n=m$ $r=m+1$, then we have

$$S_1 = \sum_{r=1}^{m+1}(-1)^{r-1} \cdot \binom{k+r-2}{r-1} \cdot \binom{m}{r-1}. \tag{111}$$

Let us apply Abel's lemma to transformation of (111):

$$S_1 = \sum_{r=1}^{m} A_r \cdot (b_r - b_{r+1}) + A_{m+1} \cdot b_{m+1},$$

where

$$A_r = \sum_{z=1}^{r}(-1)^{z-1}\cdot\binom{m}{z-1}, \quad A_{m+1} = \sum_{z=1}^{m+1}(-1)^{z-1}\cdot\binom{m}{z-1},$$

$$b_r = \binom{k+r-2}{r-1} = \binom{k+r-2}{k-1},$$

$$b_{r+1} = \binom{k+r-1}{r} = \binom{k+r-1}{k-1}.$$

In accordance with the identity

$$\binom{a}{b} - \binom{a-1}{b} \equiv \binom{a-1}{b-1}$$

we have

$$b_r - b_{r+1} = -\binom{k+r-2}{k-2} = -\binom{k+r-2}{r}. \tag{112}$$

Returning to the summation index $n$ we have $n = z-1$, $z = n+1$, for $z=1$ $n=0$, for $z=r$ $n=r-1$, then

$$A_r = \sum_{n=0}^{r-1}(-1)^n\cdot\binom{m}{n}.$$

In accordance with (106),

$$A_r = (-1)^{r-1}\cdot\binom{m-1}{r-1}. \tag{113}$$

In particulars,

$$A_{m+1} = \sum_{n=0}^{m}(-1)^n\cdot\binom{m}{n} \equiv (1-1)^m = 0. \tag{114}$$

In accordance with (112 – 114) we have

$$S_1 = -\sum_{r=1}^{m}(-1)^{r-1}\cdot\binom{k+r-2}{r}\cdot\binom{m-1}{r-1}. \tag{115}$$

Let us apply Abel's lemma to transformation of (115):

$$S_1 = -\sum_{r=1}^{m-1}A_r\cdot(b_r - b_{r+1}) - A_m\cdot b_m, \tag{116}$$

where

$$A_r = \sum_{z=1}^{r}(-1)^{z-1}\cdot\binom{m-1}{z-1}, \quad A_m = \sum_{z=1}^{m}(-1)^{z-1}\cdot\binom{m-1}{z-1},$$

$$b_r = \binom{k+r-2}{r} = \binom{k+r-2}{k-2},$$

$$b_{r+1} = \binom{k+r-1}{r+1} = \binom{k+r-1}{k-2}.$$

In accordance with the identity

$$\binom{a}{b} - \binom{a-1}{b} \equiv \binom{a-1}{b-1}$$

we have

$$b_r - b_{r+1} = -\binom{k+r-2}{k-3} = -\binom{k+r-2}{r+1}. \tag{117}$$

Returning to the summation index $n$ in equation for $A_r$ of (116) we have $n = z - 1$, $z = n + 1$, for $z = 1$ $n = 0$, for $z = r$ $n = r - 1$.

Then

$$A_r = \sum_{n=0}^{r-1} (-1)^n \cdot \binom{m-1}{n}.$$

In accordance with (106),

$$A_r = (-1)^{r-1} \cdot \binom{m-2}{r-1}. \tag{118}$$

In particulars,

$$A_m = \sum_{n=0}^{m-1} (-1)^n \cdot \binom{m-1}{n} = (-1)^{m-1} \cdot \binom{m-2}{m-1} \equiv (1-1)^{m-1} = 0. \tag{119}$$

In accordance with (117 – 119) we have

$$S_1 = \sum_{r=1}^{m-1} (-1)^{r-1} \cdot \binom{k+r-2}{r+1} \cdot \binom{m-2}{r-1}. \tag{120}$$

(111), (115) and (120) can be represented as unitary formula

$$S_1 = (-1)^a \cdot \sum_{r=1}^{m+1-a} (-1)^{r-1} \cdot \binom{k+r-2}{r-1+a} \cdot \binom{m-a}{r-1}, \tag{121}$$

correct for $a = 0,...,2$.

Let us prove that (121) is correct for any integer $a$ sufficient for existence of the sum of (121) right part.

Proposition 8:

*For integers m, k, a, $k \geq 1$, $a \geq 0$, $m \geq 0$ we have*

$$\sum_{n=0}^{m} (-1)^n \cdot \binom{m}{n} \cdot \binom{k+n-1}{n} = (-1)^a \cdot \sum_{r=1}^{m+1-a} (-1)^{r-1} \cdot \binom{k+r-2}{r-1+a} \cdot \binom{m-a}{r-1}. \tag{122}$$

Proof.

Obviously for $a = 0,...,2$ Proposition 8 is correct.

Let us use the induction hypothesis. Let Proposition 8 be correct for some integer $a > 2$ by assumption. If Proposition 8 is correct for $(a+1)$ then next equation should be correct

$$S_1 = (-1)^{a+1} \cdot \sum_{r=1}^{m-a} (-1)^{r-1} \cdot \binom{k+r-2}{r+a} \cdot \binom{m-a-1}{r-1}. \tag{123}$$

As Proposition 8 is correct for a by assumption let us transform (122) using Abel's lemma:

$$S_1 = (-1)^a \cdot \left\{ \sum_{r=1}^{m-a} A_r \cdot (b_r - b_{r+1}) + A_{m-a-1} \cdot b_{m-a-1} \right\}, \qquad (124)$$

$$A_r = \sum_{z=1}^{r} (-1)^{z-1} \cdot \binom{m-a}{z-1}, \quad A_{m-a+1} = \sum_{z=1}^{m-a+1} (-1)^{z-1} \cdot \binom{m-a}{z-1},$$

$$b_r = \binom{k+r-2}{r+a-1} = \binom{k+r-2}{k-a-1},$$

$$b_{r+1} = \binom{k+r-1}{r+a} = \binom{k+r-1}{k-a-1}.$$

In accordance with the identity

$$\binom{a}{b} - \binom{a-1}{b} \equiv \binom{a-1}{b-1},$$

we have

$$b_r - b_{r+1} = -\binom{k+r-2}{k-a-2} = -\binom{k+r-2}{r+a}. \qquad (125)$$

Returning to summation index $n$ in equation for $A_r$ of (125) we have $n = z-1$, $z = n+1$, for $z=1$ $n=0$, for $z=r$ $n=r-1$, then

$$A_r = \sum_{n=0}^{r-1} (-1)^n \cdot \binom{m-a}{n}.$$

In accordance with (106), we have

$$A_r = (-1)^r \cdot \binom{m-a-1}{r-1}. \qquad (126)$$

In particulars,

$$A_{m-a+1} = \sum_{n=0}^{m-a} (-1)^n \cdot \binom{m-a}{n} = (-1)^{m-a} \cdot \binom{m-a-1}{m-a} = (1-1)^{m-a} = 0. \qquad (127)$$

In accordance with (125 – 127), we have

$$S_1 = (-1)^{a+1} \cdot \sum_{r=1}^{m-a} (-1)^{r-1} \cdot \binom{m-a-1}{r-1} \cdot \binom{k+r-2}{r+a},$$

what coincides with (123).

Proposition 8 is proved.

Let us consider the values of a, for which right part of (123) can be calculated immediately without further Abel's transformations.

For $a = k - 1$ we have

$$S_1 = (-1)^{k-1} \cdot \sum_{r=1}^{m+1-k+1} (-1)^{r-1} \cdot \binom{k+r-2}{k+r-2} \cdot \binom{m-k+1}{r-1}.$$

In accordance with the identity

we have
$$S_1 = (-1)^{k-1} \cdot \sum_{r=1}^{m-k+2} (-1)^{r-1} \cdot \binom{m-k+1}{r-1}. \tag{128}$$

The sum (128) exists for $m - k + 2 \geq 1$, or $m - k \geq -1$, or
$$m \geq k - 1. \tag{129}$$

Returning to summation index $n$ in (128) we have $n = r - 1$, $r = n + 1$, for $r = 1$ $n = 0$, for $r = m - k + 2$ $n = m - k + 1$, then
$$S_1 = (-1)^{k-1} \cdot \sum_{n=0}^{m-k+1} (-1)^n \cdot \binom{m-k+1}{n} \equiv (-1)^{k-1} \cdot (1-1)^{m-k+1} = 0. \tag{130}$$

For $a = m$ we have
$$S_1 = (-1)^m \cdot \sum_{r=1}^{m-m+1} (-1)^0 \cdot \binom{k+1-2}{1-1+m} \cdot \binom{0}{0} = (-1)^m \cdot \binom{k-1}{m}. \tag{131}$$

(131) exists for any nonnegative $(k-1)$ and $m$. (130) is automatically obtained from (131) under the condition (129).

We have finally
$$\sum_{n=0}^{m} (-1)^n \cdot \binom{m}{n} \cdot \binom{k+n-1}{n} = (-1)^m \cdot \binom{k-1}{m}. \tag{132}$$

In accordance with (132) (104) can be represented as
$$S = \sum_{k=m+1}^{\infty} q'_k \cdot (-1)^m \cdot \binom{k-1}{m} = (-1)^m \cdot \sum_{k=m+1}^{\infty} q'_k \cdot \binom{k-1}{m}. \tag{133}$$

Let us introduce summation index $z = k - m - 1$, then $k = z + m + 1$, for $k = m + 1$ $z = 0$, then
$$S = (-1)^m \cdot \sum_{z=0}^{\infty} q'_{z+m+1} \cdot \binom{z+m}{m}. \tag{134}$$

Let us consider the behaviour of (134) in infinite increasing of $m$.
Let us note radius of convergence of (24) by $R'$.
Under Theorem 1 hypothesis $v(z)$ is analytic in region $H$, so
$$R' > 2 \tag{135}$$

In accordance with d'Alembert sign of convergence, we have for (24)
$$R' = \lim_{m \to \infty} \frac{q'_m}{q'_{m+1}} = \lim_{m \to \infty} \sqrt[k]{\frac{q'_k}{q'_{m+k}}}, \tag{136}$$

so in infinite increasing of $m$ (134) can be represented as
$$\lim_{m \to \infty} S = \lim_{m \to \infty} (-1)^m \cdot q_{m+1} \cdot \sum_{n=0}^{\infty} \binom{n+m}{n} \cdot \frac{1}{(R')^n}, \tag{137}$$

what can be transformed to

$$\lim_{m\to\infty} S = \lim_{m\to\infty}(-1)^m \cdot q_{m+1} \cdot \frac{1}{\left(1-\frac{1}{R'}\right)^{m+1}}. \tag{138}$$

In accordance with (135) we have

$$\frac{1}{1-\frac{1}{R'}} < R' \tag{139}$$

In accordance with Abel's theorem, power series converges absolutely inside its disk of convergence. So, in accordance with (139), (24) converges absolutely for

$$x = \frac{1}{1-\frac{1}{R'}}.$$

For a converging series $\sum_{n=0}^{\infty} a_n$ the equation exists

$$\lim_{m\to\infty} a_m = 0. \tag{140}$$

Then

$$\lim_{m\to\infty} \frac{q_{m+1}}{\left(1-\frac{1}{R'}\right)^{m+1}} = 0 \tag{141}$$

and

$$\lim_{m\to\infty} S = 0. \tag{142}$$

Lemma 1.4 is proved.

Lemma 1.5.

In accordance with (31) and (88) we have

$$q_1' = \sum_{n=1}^{m}\left\{\binom{m}{n+1} - m\cdot\binom{m}{n}\right\} \cdot (c_n - R_{n,m}). \tag{143}$$

In infinite increasing of dimension $m$ of approximation (5) we have

$$q_1' = \lim_{m\to\infty}\sum_{n=1}^{m}\left\{\binom{m}{n+1} - m\cdot\binom{m}{n}\right\} \cdot (c_n - R_{n,m}),$$

or

$$q_1' = \lim_{m\to\infty}\sum_{n=1}^{m}\left\{\binom{m}{n+1} - m\cdot\binom{m}{n}\right\} \cdot c_n -$$

$$- \lim_{m\to\infty}\sum_{n=1}^{m}\left\{\binom{m}{n+1} - m\cdot\binom{m}{n}\right\} \cdot R_{n,m}. \tag{144}$$

In accordance with the second equation of (88) we can represent (144) as

$$\lim_{m\to\infty} q'_{1,m} = q'_1 + \lim_{m\to\infty} \sum_{n=1}^{m}\left\{\binom{m}{n+1} - m\cdot\binom{m}{n}\right\}\cdot R_{n,m}. \qquad (145)$$

Then the proof of Lemma 1.5 can be reduced to the proof of proposition

$$\lim_{m\to\infty}\sum_{n=1}^{m}\left\{\binom{m}{n+1} - m\cdot\binom{m}{n}\right\}\cdot R_{n,m} = 0 \qquad (146)$$

under hypothesis of Theorem 1.

Let us note

$$S = \sum_{n=1}^{m}\left\{\binom{m}{n+1} - m\cdot\binom{m}{n}\right\}\cdot R_{n,m}.$$

In accordance with (94) we have

$$S = \sum_{n=1}^{m}\left\{\binom{m}{n+1} - m\cdot\binom{m}{n}\right\}\cdot(-1)^n\cdot\sum_{k=m+1}^{\infty} q'_k\cdot\binom{k+n-1}{n}. \qquad (147)$$

In accordance with (96 – 103) we can change the order of (147) summation. Then we have

$$S = \sum_{k=m+1}^{\infty} q'_k \cdot \sum_{n=1}^{m}(-1)^n\cdot\left\{\binom{m}{n+1} - m\cdot\binom{m}{n}\right\}\cdot\binom{k+n-1}{n}. \qquad (148)$$

Let us consider the sum

$$S_1 = \sum_{n=1}^{m}(-1)^n\cdot\left\{\binom{m}{n+1} - m\cdot\binom{m}{n}\right\}\cdot\binom{k+n-1}{n}. \qquad (149)$$

Let us note

$$S_1 = S_2 - S_3 \qquad (150)$$

where

$$S_2 = \sum_{n=1}^{m}(-1)^n\cdot\binom{m}{n+1}\cdot\binom{k+n-1}{n},$$

$$S_3 = \sum_{n=1}^{m}(-1)^n\cdot m\cdot\binom{m}{n}\cdot\binom{k+n-1}{n} = m\cdot\sum_{n=1}^{m}(-1)^n\cdot\binom{m}{n}\cdot\binom{k+n-1}{n}.$$

Let us consider the sums of (150). For the first sum we have

$$S_2 = \sum_{n=1}^{m}(-1)^n\cdot\binom{m}{n+1}\cdot\binom{k+n-1}{n} = \sum_{n=1}^{m}(-1)^n\cdot\binom{m}{n+1}\cdot\binom{k+n-1}{k-1}.$$

In accordance with property of binomial coefficients we have

$$S_2 = \sum_{n=1}^{m-1}(-1)^n\cdot\binom{m}{n+1}\cdot\binom{k+n-1}{k-1}.$$

Let us apply Abel's summation formula to transformation of $S_2$:

$$S_2 = \sum_{n=1}^{m-2} A_n\cdot(b_n - b_{n+1}) + b_{m-1}\cdot A_{m-1}, \qquad (151)$$

where

$$A_n = \sum_{s=1}^{n}(-1)^s \cdot \binom{m}{s+1}, \quad A_{m-1} = \sum_{s=1}^{m-1}(-1)^s \cdot \binom{m}{s+1},$$

$$b_n = \binom{k+n-1}{k-1}, \quad b_{m-1} = \binom{k+m-2}{k-1},$$

$$b_{n+1} = \binom{k+n}{k-1}.$$

In accordance with the identity

$$\binom{a}{b} - \binom{a-1}{b} \equiv \binom{a-1}{b-1}$$

we have

$$b_n - b_{n+1} = -\binom{k+n-1}{k-2}. \tag{152}$$

Let us consider the transformation of $A_n$.

Let us introduce summation index $z = s+1$, then $s = z-1$, for $s=1$ $z=2$, for $s=n$ $z=n+1$, then we have

$$A_n = \sum_{z=2}^{n+1}(-1)^{z-1} \cdot \binom{m}{z}. \tag{153}$$

(153) can be transformed to

$$A_n = -\sum_{z=2}^{n+1}(-1)^z \cdot \binom{m}{z},$$

or

$$A_n = -\left\{\sum_{z=0}^{n+1}(-1)^z \cdot \binom{m}{z} - \sum_{z=0}^{1}(-1)^z \cdot \binom{m}{z}\right\}. \tag{154}$$

In accordance with (106) we have

$$A_n = -\left\{(-1)^{n+1} \cdot \binom{m-1}{n+1} - (-1)^1 \cdot \binom{m-1}{1}\right\},$$

or

$$A_n = (-1)^n \cdot \binom{m-1}{n+1} - (m-1). \tag{155}$$

In particulars,

$$A_{m-1} = (-1)^{m-1} \cdot \binom{m-1}{m} - (m-1).$$

In accordance with the property of binomial coefficients,

$$A_{m-1} = -(m-1). \tag{156}$$

In accordance with (152 – 156) we can transform (151) to

$$S_2 = -\sum_{n=1}^{m-2}\left\{(-1)^n \cdot \binom{m-1}{n+1} - (m-1)\right\} \cdot \binom{k+n-1}{k-2} - (m-1) \cdot \binom{k+m-2}{k-1},$$

or
$$S_2 = -\sum_{n=1}^{m-2}(-1)^n \cdot \binom{m-1}{n+1} \cdot \binom{k+n-1}{k-2} + \sum_{n=1}^{m-2}(m-1) \cdot \binom{k+n-1}{k-2} - $$
$$- (m-1) \cdot \binom{k+m-2}{k-1}. \tag{157}$$

Let us consider the second additive of (157).

Let us introduce summation index $z = n+1$, then $n = z-1$, for $n=1$ $z=2$, for $n = m-2$ $z = m-1$, then we have

$$(m-1) \cdot \sum_{n=1}^{m-2} \binom{k+n-1}{k-2} = (m-1) \cdot \sum_{z=2}^{m-1} \binom{k+z-2}{k-2}. \tag{158}$$

The sum of (158) right part can be represented as

$$\sum_{z=2}^{m-1} \binom{k+z-2}{k-2} = \sum_{z=0}^{m-1} \binom{k+z-2}{k-2} - \sum_{z=0}^{1} \binom{k+z-2}{k-2}. \tag{159}$$

In accordance with [16] there exist next identities for the additives of (159) right part

$$\sum_{z=0}^{m-1} \binom{k+z-2}{k-2} = \binom{k-2+m-1+1}{k-2+1} = \binom{k+m-2}{k-1}, \tag{160}$$

$$\sum_{z=0}^{1} \binom{k+z-2}{k-2} = \binom{k-2+1-1}{k-2+1} = \binom{k}{k-1} \equiv k. \tag{161}$$

In accordance with (159 – 160) (158) can be transformed to

$$(m-1) \cdot \sum_{n=1}^{m-2} \binom{k+n-1}{k-2} = (m-1) \cdot \left\{ \binom{k+m-2}{k-1} - k \right\},$$

so (157) to

$$S_2 = -\sum_{n=1}^{m-2}(-1)^n \cdot \binom{m-1}{n+1} \cdot \binom{k+n-1}{k-2} + (m-1) \cdot \left\{ \binom{k+m-2}{k-1} - \right.$$
$$\left. - k - \binom{k+m-2}{k-1} \right\},$$

or

$$S_2 = -\sum_{n=1}^{m-2}(-1)^n \cdot \binom{m-1}{n+1} \cdot \binom{k+n-1}{k-2} + (m-1) \cdot k. \tag{162}$$

Let us apply Abel's lemma to transformation of (162):

$$S_2 = -\left\{ \sum_{n=1}^{m-3} A_n \cdot (b_n - b_{n+1}) + A_{m-2} \cdot b_{m-2} \right\} - (m-1) \cdot k, \tag{163}$$

where

$$A_n = \sum_{s=1}^{n}(-1)^s \cdot \binom{m-1}{s+1}, \quad A_{m-2} = \sum_{s=1}^{m-2}(-1)^s \cdot \binom{m-1}{s+1},$$

$$b_n = \binom{k+n-1}{k-2}, \qquad b_{m-2} = \binom{k+m-3}{k-2},$$

$$b_{n+1} = \binom{k+n}{k-2}.$$

In accordance with the identity

$$\binom{a}{b} - \binom{a-1}{b} \equiv \binom{a-1}{b-1}$$

we have

$$b_n - b_{n+1} = -\binom{k+n-1}{k-3}. \tag{164}$$

In accordance with (151 – 155) we have

$$\sum_{s=1}^{n}(-1)^s \cdot \binom{m}{s+1} = (-1)^n \cdot \binom{m-1}{n+1} - (m-1). \tag{165}$$

In accordance with (165),

$$A_n = (-1)^n \cdot \binom{m-2}{n+1} - (m-2). \tag{166}$$

In particulars,

$$A_{m-2} = (-1)^{m-2} \cdot \binom{m-2}{m-1} - (m-2).$$

In accordance with the property of binomial coefficients we have
$$A_{m-2} = -(m-2). \tag{167}$$

In accordance with (164 – 167), (163) can be transformed to

$$S_2 = \sum_{n=1}^{m-3}\left\{(-1)^n \cdot \binom{m-2}{n+1} - (m-2)\right\} \cdot \binom{k+n-1}{k-3} +$$

$$+ (m-2) \cdot \binom{k+m-3}{k-2} - (m-1) \cdot k,$$

or

$$S_2 = \sum_{n=1}^{m-3}(-1)^n \cdot \binom{m-2}{n+1}\binom{k+n-1}{k-3} - (m-2) \cdot \sum_{n=1}^{m-3}\binom{k+n-1}{k-3} +$$

$$+ (m-2) \cdot \binom{k+m-3}{k-2} - (m-1) \cdot k. \tag{168}$$

Let us consider the second additive of (168).

Let us introduce summation index $z = n+2$, then $n = z-2$, for $n=1$ $z=3$, for $n = m-3$ $z = m-1$, then we have

$$(m-2) \cdot \sum_{n=1}^{m-3}\binom{k+n-1}{k-3} = (m-2) \cdot \sum_{z=3}^{m-1}\binom{k+z-3}{k-3}. \tag{169}$$

The sum of (169) right part can be represented as

$$\sum_{z=3}^{m-1}\binom{k+z-3}{k-3}=\sum_{z=0}^{m-1}\binom{k+z-3}{k-3}-\sum_{z=0}^{2}\binom{k+z-3}{k-3}. \tag{170}$$

In accordance with (16), the identities for additives of (170) right part are

$$\sum_{z=0}^{m-1}\binom{k+z-3}{k-3}=\binom{k-3+m-1+1}{k-3+1}=\binom{k+m-3}{k-2}, \tag{171}$$

$$\sum_{z=0}^{2}\binom{k+z-3}{k-3}=\binom{k-3+2+1}{k-3+1}=\binom{k}{k-2}=\binom{k}{2}. \tag{172}$$

In accordance with (171–172), (169) can be represented as

$$(m-1)\cdot\sum_{n=1}^{m-3}\binom{k+n-1}{k-3}=(m-2)\cdot\left\{\binom{k+m-3}{k-2}-\binom{k}{2}\right\},$$

then we have next identity for (168)

$$S_2=\sum_{n=1}^{m-3}(-1)^n\cdot\binom{m-2}{n+1}\binom{k+n-1}{k-3}-$$

$$-(m-2)\cdot\left\{\binom{k+m-3}{k-2}-\binom{k}{2}-\binom{k+m-3}{k-2}\right\}-(m-1)\cdot k,$$

or

$$S_2=\sum_{n=1}^{m-3}(-1)^n\cdot\binom{m-2}{n+1}\binom{k+n-1}{k-3}-(m-2)\cdot\binom{k}{2}-(m-1)\cdot k. \tag{173}$$

(162) and (173) can be represented as unitary formula

$$S_2=(-1)^a\cdot\sum_{n=1}^{m-1-a}(-1)^n\cdot\binom{m-a}{n+1}\binom{k+n-1}{k-1-a}+\sum_{r=1}^{a}(-1)^r\cdot(m-r)\cdot\binom{k}{r}, \tag{174}$$

correct for $a=1,...,2$.

Let us prove that (174) is correct for any integer $a$ sufficient for existence of both sums of (174) right part.

Sufficient condition for the first sum is
$$m-1-a\geq 1,$$
or
$$a\leq m-2.$$
Sufficient condition for the second sum is
$$a\geq 1.$$

Proposition 9.

*For integers m, k, a, $m>1$, $a\in[1,m-2]$ we have*

$$\sum_{n=1}^{m-1}(-1)^n\cdot\binom{m}{n+1}\binom{k+n-1}{k-1}=(-1)^a\cdot\sum_{n=1}^{m-1-a}(-1)^n\cdot\binom{m-a}{n+1}\binom{k+n-1}{k-1-a}+$$

$$+\sum_{r=1}^{a}(-1)^r\cdot(m-r)\cdot\binom{k}{r}. \tag{175}$$

Proof.

Obviously for $a=1,\ldots,2$ Proposition 9 is correct.

Let us use induction hypothesis. Let proposition 9 be correct for some integer $a \in [3, m-2]$ by assumption. If Proposition 9 is correct for $(a+1)$ then next equation should be correct

$$\sum_{n=1}^{m-1}(-1)^n \cdot \binom{m}{n+1}\binom{k+n-1}{k-1} = (-1)^{a+1} \cdot \sum_{n=1}^{m-2-a}(-1)^n \cdot \binom{m-a-1}{n+1}\binom{k+n-1}{k-2-a} +$$
$$+ \sum_{r=1}^{a+1}(-1)^r \cdot (m-r) \cdot \binom{k}{r}. \tag{176}$$

As Proposition 9 is correct for $a$ by assumption let us transform (175) using Abel's lemma:

$$S_2 = (-1)^a \cdot \left\{ \sum_{n=1}^{m-2-a} A_n \cdot (b_n - b_{n+1}) + b_{m-1-a} \cdot A_{m-1-a} \right\} +$$
$$+ \sum_{r=1}^{a}(-1)^r \cdot (m-r) \cdot \binom{k}{r}, \tag{177}$$

where

$$A_n = \sum_{s=1}^{n}(-1)^s \cdot \binom{m-a}{s+1}, \quad A_{m-1-a} = \sum_{s=1}^{m-1-a}(-1)^s \cdot \binom{m-a}{s+1},$$

$$b_n = \binom{k+n-1}{k-1-a}, \quad b_{m-1-a} = \binom{k+m-1-a-1}{k-1-a} = \binom{k+m-2-a}{k-1-a},$$

$$b_{n-1} = \binom{k+n}{k-1-a}.$$

In accordance with the identity

$$\binom{a}{b} - \binom{a-1}{b} \equiv \binom{a-1}{b-1}$$

we have

$$b_n - b_{n+1} = -\binom{k+n-1}{k-2-a}. \tag{178}$$

Let us consider the transformation of $A_n$.

In accordance with (165) we have

$$A_n = (-1)^n \cdot \binom{m-a-1}{n+1} - (m-a-1). \tag{179}$$

In particulars,

$$A_{m-a-1} = (-1)^n \cdot \binom{m-a-1}{m-a} - (m-a-1).$$

In accordance with the property of binomial coefficients we have

$$A_{m-a-1} = -(m-a-1). \tag{180}$$

In accordance with (178 – 180), (177) can be transformed to

$$S_2 = (-1)^a \cdot \left\{ -\sum_{n=1}^{m-2-a} \left[ (-1)^n \cdot \binom{m-a-1}{n+1} - (m-a-1) \right] \cdot \binom{k+n-1}{k-2-a} - \right.$$

$$\left. - (m-a-1) \cdot \binom{k+m-a-2}{k-a-1} \right\} + \sum_{r=1}^{a} (-1)^r \cdot (m-r) \cdot \binom{k}{r},$$

or

$$S_2 = (-1)^{a+1} \cdot \sum_{n=1}^{m-2-a} (-1)^n \cdot \binom{m-a-1}{n+1} \cdot \binom{k+n-1}{k-a-2} -$$

$$- (-1)^{a+1} \cdot (m-a-1) \cdot \sum_{n=1}^{m-2-a} \binom{k+n-1}{k-2-a} + (-1)^{a+1} \cdot (m-a-1) \cdot \binom{k+m-a-2}{k-a-1} +$$

$$+ \sum_{r=1}^{a} (-1)^r \cdot (m-r) \cdot \binom{k}{r}. \tag{181}$$

Let us consider the second sum of (181).

Let us introduce summation index $z = n + a + 1$, then $n = z - a - 1$, for $n = 1$ $z = a + 2$, for $n = m - 2 - a$ $z = m - 1$, then we have

$$\sum_{n=1}^{m-2-a} \binom{k+n-1}{k-2-a} = \sum_{z=a+2}^{m-1} \binom{k+z-2-a}{k-2-a}. \tag{182}$$

The second additive of (182) can be transformed to

$$\sum_{z=a+2}^{m-1} \binom{k+z-2-a}{k-2-a} = \sum_{z=0}^{m-1} \binom{k+z-2-a}{k-2-a} - \sum_{z=0}^{a+1} \binom{k+z-2-a}{k-2-a}. \tag{183}$$

In accordance with (16) the identities of (183) right part additives are

$$\sum_{z=0}^{m-1} \binom{k+z-a-2}{k-a-2} = \binom{k+m-a-2}{k-a-1}, \tag{184}$$

$$\sum_{z=0}^{a+1} \binom{k+z-a-2}{k-a-2} = \binom{k}{k-a-1} \equiv \binom{k}{a+1}. \tag{185}$$

In accordance with (183 – 185), (182) can be transformed to

$$\sum_{n=1}^{m-a-2} \binom{k+n-1}{k-a-2} = \binom{k+m-a-2}{k-a-1} - \binom{k}{a+1},$$

and (181) to

$$S_2 = (-1)^{a+1} \cdot \sum_{n=1}^{m-2-a} (-1)^n \cdot \binom{m-a-1}{n+1} \cdot \binom{k+n-1}{k-a-2} -$$

$$- (-1)^{a+1} \cdot (m-a-1) \cdot \left\{ \binom{k+m-a-2}{k-a-1} - \binom{k}{a+1} - \binom{k+m-a-2}{k-a-1} \right\} +$$

$$+ \sum_{r=1}^{a} (-1)^r \cdot (m-r) \cdot \binom{k}{r},$$

or

$$S_2 = (-1)^{a+1} \cdot \sum_{n=1}^{m-2-a} (-1)^n \cdot \binom{m-a-1}{n+1} \cdot \binom{k+n-1}{k-2-a} + \sum_{r=1}^{a+1} (-1)^r \cdot (m-r) \cdot \binom{k}{r},$$

what coincides with (176).

Proposition 9 is proved.

Let $a$ in (175) be equal $(a = m - 2)$, then first additive of (175) right part can be transformed to

$$(-1)^a \cdot \sum_{n=1}^{m-1-a} (-1)^n \cdot \binom{m-a}{n+1} \cdot \binom{k+n-1}{k-1-a} =$$

$$= (-1)^{m-2} \cdot \sum_{n=1}^{1} (-1)^n \cdot \binom{2}{n+1} \cdot \binom{k+n-1}{k+m+1} =$$

$$-(-1)^{m-2} \cdot \binom{k}{k-m+1} = (-1)^{m+1} \cdot \binom{k}{m-1}. \tag{186}$$

The second additive of (175) will be

$$\sum_{r=1}^{a} (-1)^r \cdot \binom{k}{r} \cdot (m-r) = \sum_{r=1}^{m-2} (-1)^r \cdot \binom{k}{r} \cdot (m-r). \tag{187}$$

For (186) we have next equation

$$(-1)^{m+1} \cdot \binom{k}{m-1} = (-1)^{m+1} \cdot \binom{k}{m-1} \cdot (m - (m-1)),$$

then (175) can be transformed to

$$\sum_{n=1}^{m-1} (-1)^n \cdot \binom{m}{n+1} \cdot \binom{k+n-1}{k-1} = \sum_{r=1}^{m-1} (-1)^r \cdot (m-r) \cdot \binom{k}{r}, \tag{188}$$

or

$$S_2 = S_4 - S_5, \tag{189}$$

where

$$S_4 = m \cdot \sum_{r=1}^{m-1} (-1)^r \cdot \binom{k}{r},$$

$$S_5 = \sum_{r=1}^{m-1} (-1)^r \cdot r \cdot \binom{k}{r}.$$

Let us consider the sums of (189).
We have

$$\sum_{r=1}^{m-1} (-1)^r \cdot \binom{k}{r} = \sum_{r=0}^{m-1} (-1)^r \cdot \binom{k}{r} - \sum_{r=0}^{0} (-1)^r \cdot \binom{k}{r}. \tag{190}$$

The second additive of (190) is

$$\sum_{r=0}^{0} (-1)^r \cdot \binom{k}{r} = 1.$$

In accordance with (106) we have for the first additive of (190)

$$\sum_{r=0}^{m-1} (-1)^r \cdot \binom{k}{r} = (-1)^{m-1} \cdot \binom{k-1}{m-1}. \tag{191}$$

In accordance with (194 – 191), we have
$$S_4 = m \cdot \left\{(-1)^{m-1} \cdot \binom{k-1}{m-1} - 1\right\}. \qquad (192)$$

Let us use Abel's lemma to transformation of $S_5$:
$$\sum_{r=1}^{m-1}(-1)^r \cdot r \cdot \binom{k}{r} = \sum_{r=1}^{m-2} A_r \cdot (b_r - b_{r+1}) + b_{m-1} \cdot A_{m-1}, \qquad (193)$$

where
$$A_r = \sum_{s=1}^{r}(-1)^s \cdot \binom{k}{s}, \qquad A_{m-1} = \sum_{s=1}^{m-1}(-1)^s \cdot \binom{k}{s},$$
$$b_r = r, \qquad b_{m+1} = m+1,$$
$$b_{r+1} = r+1.$$

We have
$$b_r - b_{r+1} = -1. \qquad (194)$$

After replacement in (190) $r$ by $s$ we have the equation for $A_{m-1}$ in (193). In accordance with (192),
$$A_{m-1} = (-1)^{m-1} \cdot \binom{k-1}{m-1} - 1. \qquad (195)$$

After replacement in (195) $(m-1)$ by $r$ we have
$$A_r = (-1)^r \cdot \binom{k-1}{r} - 1. \qquad (196)$$

In accordance with (194 – 196), we can transform (193) to
$$S_5 = -\sum_{r=1}^{m-2}\left\{(-1)^r \cdot \binom{k-1}{r} - 1\right\} + (m-1) \cdot \left\{(-1)^{m-1} \cdot \binom{k-1}{m-1} - 1\right\},$$

or
$$S_5 = -\sum_{r=1}^{m-2}(-1)^r \cdot \binom{k-1}{r} + \sum_{r=1}^{m-2} 1 + (m-1) \cdot (-1)^{m-1} \cdot \binom{k-1}{m-1} - m + 1. \qquad (197)$$

In accordance with (106) we have next equation for the first additive of (197):
$$-\sum_{r=1}^{m-2}(-1)^r \cdot \binom{k-1}{r} = -\left\{(-1)^{m-2} \cdot \binom{k-2}{m-2} - 1\right\} = 1 - (-1)^{m-2} \cdot \binom{k-2}{m-2}. \qquad (198)$$

For the second additive we have
$$\sum_{r=1}^{m-2} 1 \equiv m - 2. \qquad (199)$$

In accordance with (198 – 199), (197) can be transformed to
$$S_5 = 1 - (-1)^{m-2} \cdot \binom{k-2}{m-2} + m - 2 + m \cdot (-1)^{m-1} \cdot \binom{k-1}{m-1} -$$
$$- (-1)^{m-1} \cdot \binom{k-1}{m-1} - m - 1.$$

In accordance with the identity

$$\binom{a}{b} - \binom{a-1}{b-1} \equiv \binom{a-1}{b}$$

we have finally for $S_5$:

$$S_5 = -(-1)^{m-1} \cdot \binom{k-2}{m-1} + m \cdot (-1)^{m-1} \cdot \binom{k-1}{m-1}. \tag{200}$$

In accordance with (199) and (200) we can represent (189) as

$$S_2 = m \cdot (-1)^{m-1} \cdot \binom{k-1}{m-1} - m + (-1)^{m-1} \cdot \binom{k-2}{m-1} - m \cdot (-1)^{m-1} \cdot \binom{k-1}{m-1},$$

or

$$S_2 = (-1)^{m-1} \cdot \binom{k-2}{m-1} - m. \tag{201}$$

Let us consider the sum $S_3$.

We have

$$S_3 = m \cdot \sum_{n=1}^{m} (-1)^n \cdot \binom{m}{n} \cdot \binom{k+n-1}{n} = m \cdot \sum_{n=1}^{m} (-1)^n \cdot \binom{m}{n} \cdot \binom{k+n-1}{k-1},$$

or

$$S_3 = m \cdot \left\{ \sum_{n=0}^{m} (-1)^n \cdot \binom{m}{n} \cdot \binom{k+n-1}{k-1} - \sum_{m=0}^{0} (-1)^n \cdot \binom{m}{n} \cdot \binom{k+n-1}{k-1} \right\}. \tag{202}$$

Let us consider the additives of the brackets of (202).

In accordance with (132),

$$\sum_{n=0}^{m} (-1)^n \cdot \binom{m}{n} \cdot \binom{k+n-1}{k-1} = (-1)^m \cdot \binom{k-1}{m}. \tag{203}$$

We have

$$\sum_{n=0}^{0} (-1)^n \cdot \binom{m}{n} \cdot \binom{k+n-1}{k-1} = 1. \tag{204}$$

In accordance with (203 – 204) we can transform (202) to

$$S_3 = m \cdot (-1)^m \cdot \binom{k-1}{m} - m. \tag{205}$$

In accordance with (201) and (205) we can transform (150) to

$$S_1 = (-1)^{m-1} \cdot \binom{k-2}{m-1} - m \cdot (-1)^m \cdot \binom{k-1}{m},$$

or

$$\sum_{n=1}^{m} \left\{ \binom{m}{n+1} - m \cdot \binom{m}{n} \right\} \cdot (-1)^n \cdot \binom{k+n-1}{n} = (-1)^{m-1} \cdot \left\{ m \cdot \binom{k-1}{m} + \binom{k-2}{m-1} \right\}. \tag{206}$$

In accordance with (206) we can transform (148) to

$$S = (-1)^{m-1} \cdot \sum_{k=m+1}^{\infty} q'_k \cdot \left\{ m \cdot \binom{k-1}{m} + \binom{k-2}{m-1} \right\}. \tag{207}$$

Let us introduce summation index $z = k - m - 1$, then $k = z + m + 1$, for $k = m+1$ $z = 0$, then we have

$$S = (-1)^{m+1} \cdot \sum_{z=0}^{\infty} q'_{z+m+1} \cdot \left\{ m \cdot \binom{z+m}{m} + \binom{z+m-1}{m-1} \right\}. \tag{208}$$

Let us consider the behaviour of (208) in unlimited increasing of $m$. Let us note radius of convergence of series (24) as $R'$. In accordance with (135 – 136) we can represent (208) in infinite increasing or $m$ as

$$\lim_{m \to \infty} S = \lim_{m \to \infty} (-1)^m \cdot q'_{m+1} \cdot \sum_{n=0}^{\infty} \left\{ m \cdot \binom{n+m}{m} + \binom{n+m-1}{m-1} \right\} \cdot \frac{1}{R'^n}. \tag{209}$$

Let us transform (209) to

$$\lim_{m \to \infty} S = \lim_{m \to \infty} (-1)^m \cdot q'_{m+1} \cdot \left\{ \frac{m}{\left(1 - \frac{1}{R'}\right)^{m+1}} + \frac{1}{\left(1 - \frac{1}{R'}\right)^m} \right\}. \tag{210}$$

Let us consider the behaviour of absolute value of (210). We have

$$\lim_{m \to \infty} |S| = \lim_{m \to \infty} \frac{q'_{m+1} \cdot m}{\left(1 - \frac{1}{R'}\right)^{m+1}} + \lim_{m \to \infty} \frac{q'_{m+1}}{\left(1 - \frac{1}{R'}\right)^m}, \tag{211}$$

or

$$C = \lim_{m \to \infty} |S| = \lim_{m \to \infty} C(m),$$

$$C_1 = \lim_{m \to \infty} \frac{q'_{m+1} \cdot m}{\left(1 - \frac{1}{R'}\right)^{m+1}} = \lim_{m \to \infty} C_1(m),$$

$$C_2 = \lim_{m \to \infty} \frac{q'_{m+1}}{\left(1 - \frac{1}{R'}\right)^m} = \lim_{m \to \infty} C_2(m).$$

Let us represent $C_1$ as

$$C_1 = \lim_{m \to \infty} \left\{ \left(1 - \frac{1}{R'}\right)^2 \cdot \frac{1}{R'} \cdot \frac{q'_m \cdot m}{\left(1 - \frac{1}{R'}\right)^{m-1}} \right\}. \tag{212}$$

In accordance with (135), (139) and Abel's theorem, for

$$x = \frac{1}{1 - \frac{1}{R'}}$$

(24) and all its derivatives converges absolutely. For converging series $\sum_{n=0}^{\infty} a_n$ (140) is correct.

First derivative of (24) can be represented as

$$\frac{d}{dx} \cdot \sum_{n=0}^{\infty} q'_n \cdot x^n = \sum_{n=1}^{\infty} q'_n \cdot n \cdot x^{n-1}. \tag{213}$$

Then

$$\lim_{m \to \infty} \frac{q'_m \cdot m}{\left(1 - \frac{1}{R'}\right)^{m-1}} = 0, \tag{214}$$

and, in accordance with (135),

$$\lim_{m \to \infty} C_1(m) = 0. \tag{215}$$

Similarly to deduction of (142) we can prove that

$$\lim_{m \to \infty} C_2(m) = 0. \tag{216}$$

Theorem 1 is proved.

## Conclusion

Theorem 1 contains sufficient conditions of existence of asymptotic evaluation (3) of a function, defined by a power series different from conditions of [16]. Nevertheless there may exist other sufficient conditions. The problem of coefficients of higher orders will be considered in the next paper.